\documentclass[10 pt]{amsart}
\usepackage[english]{babel}
\usepackage{graphicx}
\usepackage{amsmath}
\usepackage{amsthm}
\usepackage{amssymb}
\usepackage{textcase}
\usepackage{stmaryrd}
\usepackage{relsize}
\usepackage[utf8]{inputenc}
\usepackage{hyperref}
\usepackage{url}
\usepackage[right=1 in, left = 1 in]{geometry}

\newtheorem{proposition}{Proposition}

\newtheorem{lemma}[proposition]{Lemma}
\newtheorem{corollary1}{Corollary}[proposition]
\newtheorem{corollary}[corollary1]{Corollary}

\theoremstyle{remark}

\title{On Some Doubly Logarithmic Integrals}

\author{Duc Van Khanh Tran}\thanks{Contact D. V. K. Tran: duc.tranvk@utexas.edu; Department of Mathematics, University of Texas at Austin, USA}

\begin{document}
\subjclass{33B15, 11M06}
\keywords{Gradshteyn and Ryzhik, logarithmic integrals, gamma function}
\maketitle

\begin{abstract}
\normalsize
\noindent
There have been many works on proving the integrals in the table of integrals compiled by Gradshteyn and Ryzhik, and in this paper we prove some doubly logarithmic integral identities in the Gradshteyn and Ryzhik table.
\end{abstract}

\section{Introduction}
There has been extensive research done on the integrals in the well-known Gradshteyn and Ryzhik table \cite{GradshteynRyzhik}. Moll and his collaborators have written a series of papers on several classes of integrals in the Gradshteyn and Ryzhik table (see, for example, \cite{Moll30}, \cite{Moll31}, and \cite{Moll32}). In \cite{boyadzhiev}, Boyadzhiev used Poisson's integration formula to evaluate some integrals appearing in the Gradshteyn and Ryzhik table. Reynolds and Stauffer evaluated some integrals involving hyperbolic functions in the Gradshteyn and Ryzhik table in \cite{robertallan1} and \cite{robertallan2}, and Coffey evaluated some integrals involving hyperbolic and trigonometric functions in the Gradshteyn and Ryzhik table in \cite{coffey}.

In \cite{Vardipaper}, Vardi proved the integral identity
$$\int_{\pi /4}^{\pi/2} \ln{\ln{\tan{x}}}\,dx = \frac{\pi}{2}\ln{\left(\frac{\Gamma\left(\frac{3}{4}\right)}{\Gamma\left(\frac{1}{4}\right)}\sqrt{2\pi}\right)},$$
which is the identity 4.229.7 in the Gradshteyn and Ryzhik table (see \cite[4.229.7]{GradshteynRyzhik}), using an analytic number theoretical method. Here, $\Gamma(\cdot)$ denotes the gamma function (see \cite[sec. 5.2]{DLMF}). The Gradshteyn and Ryzhik table also lists the identity
$$\int_0^1 \ln{(-\ln{x})}\frac{1}{1+x^2}\,dx = \frac{\pi}{2}\ln{\frac{\sqrt{2\pi}\Gamma\left(\frac{3}{4}\right)}{\Gamma\left(\frac{1}{4}\right)}}$$
(see \cite[4.325.4]{GradshteynRyzhik}). Note that these identities are equivalent since substituting $\cot{x} \to x$ gives
$$\int_{\pi /4}^{\pi/2} \ln{\ln{\tan{x}}}\,dx = \int_0^1 \ln{(-\ln{x})}\frac{1}{1+x^2}\,dx.$$
In \cite{Vardipaper}, Vardi also described how to evaluate two other similar identities,
$$\int_0^1 \ln{(-\ln{x})} \frac{1}{1+x+x^2}\,dx = \frac{\pi}{\sqrt{3}}\ln{\frac{\sqrt[3]{2\pi}\Gamma\left(\frac{2}{3}\right)}{\Gamma\left(\frac{1}{3}\right)}}$$
and
$$\int_0^1 \ln{(-\ln{x})}\frac{1}{1-x+x^2}\,dx = \frac{2\pi}{\sqrt{3}}\left[\frac{5}{6}\ln{2\pi} - \ln{\Gamma\left(\frac{1}{6}\right)}\right],$$
which are the identities 4.325.5 and 4.325.6 in the Gradshteyn and Ryzhik table (see \cite[4.325.5-6]{GradshteynRyzhik}), by using the same number theoretical method. In Section 2, we prove these integral identities analytically without involving number theoretical methods. Additionally, we prove a few other doubly logarithmic integral identities listed in the Gradshteyn and Ryzhik table:
$$\int_0^1 \ln{(-\ln{x})} \frac{1}{1+x} \,dx= -\frac{1}{2}(\ln{2})^2,$$
$$\int_0^1 \ln{(-\ln{x})}\frac{1}{x+e^{i\lambda}} \, dx= \sum_{k=1}^{\infty} \frac{(-1)^k}{k}e^{-ik\lambda}(\gamma + \ln{k}),$$
$$\int_0^1 \ln{(-\ln{x})}\frac{1}{(1+x)^2}\,dx = \frac{1}{2}\left(\ln{\frac{\pi}{2}} - \gamma\right),$$
$$\int_0^1 \ln{(-\ln{x})}\frac{1}{1+2x\cos{t} + x^2}\,dx = \frac{\pi}{2\sin{t}} \ln{\frac{(2\pi)^{t/\pi}\Gamma\left(\frac{1}{2}+\frac{t}{2\pi}\right)}{\Gamma\left(\frac{1}{2}-\frac{t}{2\pi}\right)}},$$
$$\int_0^1 \ln{(-\ln{x})} \frac{1}{\left(1+x^2\right)\sqrt{-\ln{x}}}\,dx = \sqrt{\pi}\sum_{k=0}^{\infty} \frac{(-1)^{k+1}}{\sqrt{2k+1}}[\ln{(2k+1) + 2\ln{2} + \gamma}]$$
(see \cite[4.325.1-3, 7, and 10]{GradshteynRyzhik}). Here, $\gamma$ denotes the Euler-Mascheroni constant defined as
$$\gamma := \lim_{n\to \infty} \left(\sum_{k=1}^n \frac{1}{k} - \ln{n}\right)$$
(see \cite[sec. 1.5]{mathconstant}). Note that the Euler-Mascheroni constant is denoted as $C$ in the Gradshteyn and Ryzhik table \cite{GradshteynRyzhik}. The integrals evaluated in this paper can also be found in Tables 147 and 148 of \textit{Nouvelles tables d’intégrales définies} \cite{BI} and Table 325 of \textit{Integraltafel Zweiter Teil: Bestimmte Integrale} \cite{GW}.

In Section 3, some possible directions for further research related to the class of doubly logarithmic integrals discussed in this paper are proposed.

\section{Evaluation of Integrals}
Before proving the integral identities 4.325.1-7, we need the following lemma.
Note that this lemma is also listed as the identity 4.325.8 in the Gradshteyn and Ryzhik table \cite{GradshteynRyzhik}.
\begin{lemma} For $\text{Re}(\mu) > 0$,
$$\int_0^1 \ln{(-\ln{x})} \,x^{\mu -1}\,dx = -\frac{1}{\mu}\left(\gamma + \ln{\mu}\right).$$
\end{lemma}
\begin{proof}
    Substituting $-\mu \ln{x} \to x$,
    \begin{align*}
        \int_0^1 \ln{(-\ln{x})} \,x^{\mu -1}\,dx &= \frac{1}{\mu}\int_0^{\infty} e^{-x}\ln{\frac{x}{\mu}}\, dx \\
        &= \frac{1}{\mu}\int_0^{\infty} e^{-x}\ln{x}\, dx - \frac{\ln{\mu}}{\mu} \int_0^{\infty} e^{-x}\,dx \\
        &= \frac{1}{\mu}\int_0^{\infty} e^{-x}\ln{x}\, dx - \frac{\ln{\mu}}{\mu}.
    \end{align*}
    It is well-known that
    $$\int_0^{\infty} e^{-x}\ln{x}\, dx = -\gamma$$
    (see \cite[sec. 9.3]{irrint} and \cite[sec. 4.3]{intsumfiesta}). Therefore,
    $$\int_0^1 \ln{(-\ln{x})} \,x^{\mu -1}\,dx = -\frac{1}{\mu}\left(\gamma + \ln{\mu}\right).$$
\end{proof}

Now we prove the identity 4.325.2 in the Gradshteyn and Ryzhik table \cite{GradshteynRyzhik}.

\begin{proposition} For $\lambda \in \mathbb{R}$,
    $$\int_0^1 \ln{(-\ln{x})}\frac{1}{x+e^{i\lambda}} \, dx= \sum_{k=1}^{\infty} \frac{(-1)^k}{k}e^{-ik\lambda}(\gamma + \ln{k}),$$
where $i = \sqrt{-1}$.
\end{proposition}
\begin{proof}
    First,
    $$\int_0^1 \ln{(-\ln{x})}\frac{1}{x+e^{i\lambda}} \, dx = e^{-i\lambda}\int_0^1 \ln{(-\ln{x})}\frac{1}{e^{-i\lambda}x+1} \, dx.$$
    Since $|e^{-i\lambda}x| = |x| \leq 1$ on $[0,1]$, $\displaystyle{\sum_{k=0}^{\infty}} (-1)^ke^{-ik\lambda}x^k$ converges to $\frac{1}{e^{-i\lambda}x+1}$ almost everywhere on $[0,1]$, and
    $$\left|\sum_{k=0}^{n} (-1)^ke^{-ik\lambda}x^k\right| \leq \sum_{k=0}^{n} \left|(-1)^ke^{-ik\lambda}x^k\right| = \sum_{k=0}^{n} x^k < \frac{1}{1-x}$$
    for all $n\in \mathbb{N}$ almost everywhere on $[0,1]$. Thus, by Lebesgue dominated convergence theorem,
    \begin{align*}
        e^{-i\lambda} \int_0^1 \ln{(-\ln{x})}\frac{1}{e^{-i\lambda}x+1} \, dx &= e^{-i\lambda} \int_0^1 \ln{(-\ln{x})} \sum_{k=0}^{\infty} (-1)^ke^{-ik\lambda}x^k\, dx \\
        &= e^{-i\lambda} \int_0^1 \ln{(-\ln{x})} \lim_{n \to \infty }\sum_{k=0}^{n} (-1)^ke^{-ik\lambda}x^k\, dx \\
        &= e^{-i\lambda} \lim_{n \to \infty}\sum_{k=0}^{n}  \int_0^1 \ln{(-\ln{x})} \, (-1)^ke^{-ik\lambda}x^k\, dx\\
        &= \sum_{k=0}^{\infty} (-1)^k e^{-i(k+1)\lambda} \int_0^1 \ln{(-\ln{x})} \, x^k \, dx\\
        &= \sum_{k=1}^{\infty} (-1)^{k-1} e^{-ik\lambda} \int_0^1 \ln{(-\ln{x})} \, x^{k-1} \, dx.
    \end{align*}
    By Lemma 1,
    $$\int_0^1 \ln{(-\ln{x})} \, x^{k-1} \, dx = -\frac{1}{k}(\gamma + \ln{k}).$$
    Therefore,
    \begin{align*}
        \int_0^1 \ln{(-\ln{x})}\frac{1}{x+e^{i\lambda}} \, dx &= \sum_{k=1}^{\infty} (-1)^{k-1} e^{-ik\lambda}  \cdot -\frac{1}{k}(\gamma + \ln{k}) \\
        &= \sum_{k=1}^{\infty} \frac{(-1)^k}{k}e^{-ik\lambda}(\gamma + \ln{k}).
    \end{align*}
\end{proof}

As a corollary, we prove the identity 4.325.1 in the Gradshteyn and Ryzhik table \cite{GradshteynRyzhik}.

\begin{corollary}
    $$\int_0^1 \ln{(-\ln{x})} \frac{1}{1+x} \,dx= -\frac{1}{2}(\ln{2})^2.$$
\end{corollary}
\begin{proof}
    Let $\lambda = 0$ in Proposition 2,
    \begin{align*}
        \int_0^1 \ln{(-\ln{x})} \frac{1}{1+x} \,dx &= \sum_{k=1}^{\infty} \frac{(-1)^k}{k} (\gamma + \ln{k}) \\
        &= \gamma\sum_{k=1}^{\infty} \frac{(-1)^k}{k} + \sum_{k=1}^{\infty} \frac{(-1)^k\ln{k}}{k}\\
        &= -\gamma\,\eta(1) + \eta'(1),
    \end{align*}
    where $\eta(\cdot)$ denotes the eta function defined as
    $$\eta(s) = \sum_{k=1}^{\infty} \frac{(-1)^{k-1}}{k^s}$$
    (see \cite[eq. 23.2.19]{FuncHandbook}). It is known that $\eta(1) = \ln{(2)}$ and that $\eta'(s)$ satisfies the relation
    \begin{equation}
        \eta'(s) = 2^{1-s}(\ln{2})\zeta(s) + \left(1-2^{1-s}\right)\zeta'(s),
    \end{equation}
    where $\zeta(\cdot)$ denotes the zeta function (see \cite{dirichleteta}). The Laurent series expansion of $\zeta(s)$ at $s=1$ is
    \begin{align}
        \zeta(s) &= \frac{1}{s-1} + \sum_{n=0}^{\infty} \frac{(-1)^n}{n!}\gamma_n (s-1)^n \\
        &= \frac{1}{s-1} + \gamma -\gamma_1(s-1) + O\left((s-1)^2\right) \nonumber,
    \end{align}
    where $\gamma_n$ are the Stieltjes constants and $\gamma_0 = \gamma$ (see \cite[eq. 25.2.4]{DLMF}). Taking the derivative,
    \begin{equation}
        \zeta'(s) = -\frac{1}{(s-1)^2} -\gamma_1 + O(s-1).
    \end{equation}
    Also, the Laurent series expansion of $2^{1-s}$ at $s=1$ is
    \begin{equation}
        2^{1-s} = 1 - (\ln{2}) (s-1) + \frac{(\ln{2})^2}{2} (s-1)^2 + O\left((s-1)^3\right).
    \end{equation}
    Substituting (2), (3), and (4) into (1) gives
    \begin{align*}
        \eta'(s) &= \left( \frac{\ln{2}}{s-1} - (\ln{2})^2 + \gamma\ln{2} + O(s-1)\right) + \left(-\frac{\ln{2}}{s-1} + \frac{(\ln{2})^2}{2} + O(s-1)\right) \\
        &= \gamma\ln{2} -\frac{(\ln{2})^2}{2} + O(s-1).
    \end{align*}
    Let $s=1$,
    $$\eta'(1) = \gamma\ln{2} -\frac{(\ln{2})^2}{2}.$$
    Therefore,
    $$\int_0^1 \ln{(-\ln{x})} \frac{1}{1+x} = -\gamma\ln{2} + \gamma\ln{2} -\frac{(\ln{2})^2}{2} = -\frac{1}{2}(\ln{2})^2.$$
\end{proof}

Next, we prove the identity 4.325.7 in the Gradshteyn and Ryzhik table \cite{GradshteynRyzhik}.

\begin{proposition} For $t \in (-\pi, \pi)$,
    $$\int_0^1 \ln{(-\ln{x})}\frac{1}{1+2x\cos{t} + x^2}\,dx = \frac{\pi}{2\sin{t}} \ln{\frac{(2\pi)^{t/\pi}\Gamma\left(\frac{1}{2}+\frac{t}{2\pi}\right)}{\Gamma\left(\frac{1}{2}-\frac{t}{2\pi}\right)}}.$$
\end{proposition}
\begin{proof}
    Since
    $$\cos{t} = \frac{e^{it} + e^{-it}}{2}$$
    where $i = \sqrt{-1}$,
    \begin{align*}
        &\int_0^1 \ln{(-\ln{x})}\frac{1}{1+2x\cos{t} + x^2}\,dx \\
        =& \int_0^1 \ln{(-\ln{x})} \frac{1}{1 + \left(e^{it} + e^{-it}\right)x + x^2}\, dx \\
        =& \int_0^1 \ln{(-\ln{x})} \frac{1}{\left(x+e^{it}\right)\left(x+e^{-it}\right)} \, dx \\
        = &\frac{1}{e^{it}-e^{-it}} \left(\int_0^1 \ln{(-\ln{x})} \frac{1}{x+e^{-it}} \,dx -  \int_0^1 \ln{(-\ln{x})} \frac{1}{x+e^{it}}\, dx\right).
    \end{align*}
    By Proposition 2,
    $$\int_0^1 \ln{(-\ln{x})}\frac{1}{x+e^{-it}} \, dx= \sum_{k=1}^{\infty} \frac{(-1)^k}{k}e^{ikt}(\gamma + \ln{k}),$$
    and
    $$\int_0^1 \ln{(-\ln{x})}\frac{1}{x+e^{it}} \, dx= \sum_{k=1}^{\infty} \frac{(-1)^k}{k}e^{-ikt}(\gamma + \ln{k}).$$
    So,
    $$\int_0^1 \ln{(-\ln{x})}\frac{1}{1+2x\cos{t} + x^2}\,dx = \frac{1}{e^{it}-e^{-it}} \sum_{k=1}^{\infty} \frac{(-1)^k}{k}\left(e^{ikt} - e^{-ikt}\right)(\gamma + \ln{k}).$$
    Since
    $$\sin{\theta} = \frac{e^{i\theta} - e^{-i\theta}}{2i},$$
    we obtain
    \begin{equation}
        \int_0^1 \ln{(-\ln{x})}\frac{1}{1+2x\cos{t} + x^2}\,dx = \frac{1}{\sin{t}}\sum_{k=1}^{\infty} \frac{(-1)^k}{k}(\gamma + \ln{k})\sin{(kt)}.
    \end{equation}
    The Fourier sine series of $\ln{\Gamma(x)}$ for $0<x<1$ is
    \begin{equation}
        \ln{\Gamma(x)} = \frac{1}{2} \ln{\frac{\pi}{\sin{(\pi x)}}} + \frac{1}{\pi}\sum_{k=1}^{\infty} \frac{(\gamma + \ln{2\pi} + \ln{k})\sin{(2k\pi x)}}{k}
    \end{equation}
    (see \cite[sec. 12.3]{modernanal}), or equivalently
    $$\ln{\Gamma(x)} = \frac{1}{2} \ln{\frac{\pi}{\sin{(\pi x)}}} + (\gamma + \ln{2\pi})\left(\frac{1}{2} - x\right) + \frac{1}{\pi}\sum_{k=1}^{\infty} \frac{\ln{k}}{k}\sin{(2k\pi x)}$$
    since
    \begin{equation}
        \frac{1}{2} - x = \frac{1}{\pi}\sum_{k=1}^{\infty} \frac{\sin{(2k\pi x)}}{k}
    \end{equation}
    for $0<x<1$ (see \cite{Kummer}). Note that
    $$0 < \frac{1}{2} + \frac{t}{2\pi} < 1$$
    when $-\pi < t <\pi$. Let $x = \frac{1}{2} + \frac{t}{2\pi}$ in (6),
    $$\ln{\Gamma\left(\frac{1}{2} + \frac{t}{2\pi}\right)} = \frac{1}{2} \ln{\frac{\pi}{\sin{\left(\frac{\pi}{2} + \frac{t}{2}\right)}}} + \frac{1}{\pi}\sum_{k=1}^{\infty} \frac{(\gamma + \ln{2\pi} + \ln{k})(-1)^k\sin{(kt)}}{k},$$
    so
    $$\sum_{k=1}^{\infty} \frac{(-1)^k}{k}(\gamma + \ln{k})\sin{(kt)} = \pi \ln{\Gamma\left(\frac{1}{2} + \frac{t}{2\pi}\right)} - \frac{\pi}{2} \ln{\frac{\pi}{\sin{\left(\frac{\pi}{2} + \frac{t}{2}\right)}}} - \ln{2\pi} \, \sum_{k=1}^{\infty} \frac{(-1)^k\sin{(kt)}}{k}.$$
    Let $x = \frac{1}{2} + \frac{t}{2\pi}$ in (7),
    $$-\frac{t}{2} = \sum_{k=1}^{\infty} \frac{(-1)^k\sin{(kt)}}{k}.$$
    Thus,
    \begin{align*}
        \sum_{k=1}^{\infty} \frac{(-1)^k}{k}(\gamma + \ln{k})\sin{(kt)}  &= \pi \ln{\Gamma\left(\frac{1}{2} + \frac{t}{2\pi}\right)} - \frac{\pi}{2} \ln{\frac{\pi}{\sin{\left(\frac{\pi}{2} + \frac{t}{2}\right)}}} + \frac{t}{2}\ln{2\pi} \\
        &= \frac{\pi}{2}\ln{\frac{(2\pi)^{t/\pi}\Gamma^2\left(\frac{1}{2} + \frac{t}{2\pi}\right)\sin{\left(\frac{\pi}{2} + \frac{t}{2}\right)}}{\pi}}.
    \end{align*}
    By the reflection formula of the gamma function,
    $$\frac{\pi}{\sin{\left(\frac{\pi}{2} + \frac{t}{2}\right)}} = \Gamma\left(\frac{1}{2} + \frac{t}{2\pi}\right)\Gamma\left(\frac{1}{2} - \frac{t}{2\pi}\right)$$
    (see \cite[sec. 5.5]{DLMF}). Therefore,
    $$\sum_{k=1}^{\infty} \frac{(-1)^k}{k}(\gamma + \ln{k})\sin{(kt)} = \frac{\pi}{2}\ln{\frac{(2\pi)^{t/\pi}\Gamma\left(\frac{1}{2}+\frac{t}{2\pi}\right)}{\Gamma\left(\frac{1}{2}-\frac{t}{2\pi}\right)}}.$$
    Substituting this into (5) gives
    $$\int_0^1 \ln{(-\ln{x})}\frac{1}{1+2x\cos{t} + x^2}\,dx = \frac{\pi}{2\sin{t}} \ln{\frac{(2\pi)^{t/\pi}\Gamma\left(\frac{1}{2}+\frac{t}{2\pi}\right)}{\Gamma\left(\frac{1}{2}-\frac{t}{2\pi}\right)}}.$$
\end{proof}

As corollaries, we prove the identities 4.325.3-6 in the Gradshteyn and Ryzhik table \cite{GradshteynRyzhik}. Note that the identity 4.325.4 is equivalent to the identity 4.229.7.

\begin{corollary}
    $$\int_0^1 \ln{(-\ln{x})}\frac{1}{(1+x)^2}\,dx = \frac{1}{2}\left(\ln{\frac{\pi}{2}} - \gamma\right).$$
\end{corollary}
\begin{proof}
    By letting $t$ approach $0$ in Proposition 3 and using L'Hôpital's rule, we get
    \begin{align*}
        \int_0^1 \ln{(-\ln{x})}\frac{1}{1+2x + x^2}\,dx &= \frac{\pi}{2} \lim_{t \to 0} \frac{1}{\cos{t}}\left[\frac{\ln{2\pi}}{\pi} + \frac{1}{2\pi}\psi\left(\frac{1}{2} + \frac{t}{2\pi}\right) + \frac{1}{2\pi} \psi\left(\frac{1}{2} - \frac{t}{2\pi}\right)\right] \\
        &= \frac{1}{2}\left[\ln{2\pi} + \psi\left(\frac{1}{2}\right)\right],
    \end{align*}
    where $\psi(\cdot)$ denotes the digamma function defined as $\psi(s) = \frac{\Gamma'(s)}{\Gamma(s)}$ (see \cite[sec. 5.2]{DLMF}). Since $\psi\left(\frac{1}{2}\right) = -\gamma - 2\ln{2}$ (see \cite[sec. 5.4]{DLMF}),
    $$\int_0^1 \ln{(-\ln{x})}\frac{1}{(1+x)^2}\,dx = \int_0^1 \ln{(-\ln{x})}\frac{1}{1+2x + x^2}\,dx = \frac{1}{2}\left(\ln{\frac{\pi}{2}} - \gamma\right).$$
\end{proof}

\begin{corollary}
    $$\int_0^1 \ln{(-\ln{x})}\frac{1}{1+x^2}\,dx = \frac{\pi}{2}\ln{\frac{\sqrt{2\pi}\Gamma\left(\frac{3}{4}\right)}{\Gamma\left(\frac{1}{4}\right)}}.$$
\end{corollary}
\begin{proof}
    Let $t = \frac{\pi}{2}$ in Proposition 3,
    $$\int_0^1 \ln{(-\ln{x})}\frac{1}{1+x^2}\,dx = \frac{\pi}{2}\ln{\frac{\sqrt{2\pi}\Gamma\left(\frac{3}{4}\right)}{\Gamma\left(\frac{1}{4}\right)}}.$$
\end{proof}

\begin{corollary}
    $$\int_0^1 \ln{(-\ln{x})} \frac{1}{1+x+x^2}\,dx = \frac{\pi}{\sqrt{3}}\ln{\frac{\sqrt[3]{2\pi}\Gamma\left(\frac{2}{3}\right)}{\Gamma\left(\frac{1}{3}\right)}}.$$
\end{corollary}
\begin{proof}
    Let $t = \frac{\pi}{3}$ in Proposition 3,
    $$\int_0^1 \ln{(-\ln{x})} \frac{1}{1+x+x^2}\,dx = \frac{\pi}{\sqrt{3}}\ln{\frac{\sqrt[3]{2\pi}\Gamma\left(\frac{2}{3}\right)}{\Gamma\left(\frac{1}{3}\right)}}.$$
\end{proof}

\begin{corollary}
    $$\int_0^1 \ln{(-\ln{x})}\frac{1}{1-x+x^2}\,dx = \frac{2\pi}{\sqrt{3}}\left[\frac{5}{6}\ln{2\pi} - \ln{\Gamma\left(\frac{1}{6}\right)}\right].$$
\end{corollary}
\begin{proof}
    Let $t = \frac{2\pi}{3}$ in Proposition 3,
    $$\int_0^1 \ln{(-\ln{x})}\frac{1}{1-x+x^2}\,dx  = \frac{\pi}{\sqrt{3}}\ln{\frac{(2\pi)^{2/3}\Gamma\left(\frac{5}{6}\right)}{\Gamma\left(\frac{1}{6}\right)}}.$$
    By the reflection formula of the gamma function,
    $$\Gamma\left(\frac{1}{6}\right)\Gamma\left(\frac{5}{6}\right) = \frac{\pi}{\sin{\frac{\pi}{6}}} = 2\pi.$$
    Therefore,
    $$\int_0^1 \ln{(-\ln{x})}\frac{1}{1-x+x^2}\,dx = \frac{\pi}{\sqrt{3}}\ln{\frac{(2\pi)^{5/3}}{\Gamma^2 \left(\frac{1}{6}\right)}} = \frac{2\pi}{\sqrt{3}}\left[\frac{5}{6}\ln{2\pi} - \ln{\Gamma\left(\frac{1}{6}\right)}\right].$$
\end{proof}
Finally, we prove a lemma and then use it to prove the identity 4.325.10 in the Gradshteyn and Ryzhik table \cite{GradshteynRyzhik}. Note that the lemma is also listed as the identity 4.325.11 in the Gradshteyn and Ryzhik table \cite{GradshteynRyzhik}. We also rewrite the infinite series in the identity 4.325.10 as a closed form involving $\beta\left(\frac{1}{2}\right)$, where $\beta(\cdot)$ denotes the Dirichlet beta function defined as
$$\beta(s) = \sum_{k=0}^{\infty} \frac{(-1)^k}{(2k+1)^s}$$
(see \cite[eq. 23.2.21]{FuncHandbook}).

\begin{lemma} For $\text{Re}(\mu) >0$,
    $$\int_0^1 \ln{(-\ln{x})} \frac{x^{\mu-1}}{\sqrt{-\ln{x}}}\,dx = - (\gamma + \ln{4\mu})\sqrt{\frac{\pi}{\mu}}.$$
\end{lemma}
\begin{proof}
    Substituting $-\mu \ln{x} \to x$ gives
    \begin{align*}
        \int_0^1 \ln{(-\ln{x})} \frac{x^{\mu-1}}{\sqrt{-\ln{x}}}\,dx &= \frac{1}{\sqrt{\mu}} \int_0^{\infty} x^{-\frac{1}{2}}e^{-x} \ln{\frac{x}{\mu}}\,dx \\
        &= \frac{1}{\sqrt{\mu}} \int_0^{\infty} x^{-\frac{1}{2}}e^{-x} \ln{x}\, dx - \frac{\ln{\mu}}{\sqrt{\mu}} \int_0^{\infty} x^{-\frac{1}{2}}e^{-x} \,dx.
    \end{align*}
    Differentiating the gamma function,
    $$\Gamma(s) = \int_0^{\infty} x^{s-1}e^{-x}\,dx,$$
    gives
    $$\int_0^{\infty} x^{s-1}e^{-x}\ln{x}\, dx = \Gamma'(s) = \psi(s)\Gamma(s).$$
    Thus,
    $$\int_0^1 \ln{(-\ln{x})} \frac{x^{\mu-1}}{\sqrt{-\ln{x}}}\,dx = \frac{1}{\sqrt{\mu}}\psi\left(\frac{1}{2}\right)\Gamma\left(\frac{1}{2}\right) - \frac{\ln{\mu}}{\sqrt{\mu}}\,\Gamma\left(\frac{1}{2}\right).$$
    Since $\Gamma\left(\frac{1}{2}\right) = \sqrt{\pi}$ and $\psi\left(\frac{1}{2}\right) = -\gamma - 2\ln{2}$ (see \cite[sec. 5.4]{DLMF}),
    $$\int_0^1 \ln{(-\ln{x})} \frac{x^{\mu-1}}{\sqrt{-\ln{x}}}\,dx = -(\gamma + \ln{4\mu})\sqrt{\frac{\pi}{\mu}}.$$
\end{proof}

\begin{proposition}
    \begin{align*}
        \int_0^1 \ln{(-\ln{x})} \frac{1}{\left(1+x^2\right)\sqrt{-\ln{x}}}\,dx &= \sqrt{\pi}\sum_{k=0}^{\infty} \frac{(-1)^{k+1}}{\sqrt{2k+1}}[\ln{(2k+1) + 2\ln{2} + \gamma}]\\
        &= -\left(\ln{\sqrt{\frac{8}{\pi}}} + \frac{\pi}{4} + \frac{\gamma}{2}\right)\sqrt{\pi}\,\beta\left(\frac{1}{2}\right).
    \end{align*}
\end{proposition}
\begin{proof}
    Since $\left|x^2\right| \leq 1$ on $[0,1]$, $\displaystyle{\sum_{k=0}^{\infty}} (-1)^kx^{2k}$ converges to $\frac{1}{1+x^2}$ almost everywhere on $[0,1]$, and
    $$\left|\sum_{k=0}^{n} (-1)^kx^{2k}\right| \leq \sum_{k=0}^{n} \left|(-1)^kx^{2k}\right| = \sum_{k=0}^{n} x^{2k} < \frac{1}{1-x^2}$$
    for all $n\in \mathbb{N}$ almost everywhere on $[0,1]$. Thus, by Lebesgue dominated convergence theorem,
    \begin{align*}
        \int_0^1 \ln{(-\ln{x})}\frac{1}{\left(1+x^2\right)\sqrt{-\ln{x}}} \, dx &= \int_0^1 \ln{(-\ln{x})} \frac{1}{\sqrt{-\ln{x}}}\sum_{k=0}^{\infty} (-1)^kx^{2k}\, dx \\
        &= \int_0^1 \ln{(-\ln{x})}  \frac{1}{\sqrt{-\ln{x}}} \lim_{n \to \infty }\sum_{k=0}^{n} (-1)^kx^{2k}\, dx \\
        &= \lim_{n \to \infty}\sum_{k=0}^{n}  \int_0^1 \ln{(-\ln{x})} \frac{(-1)^kx^{2k}}{\sqrt{-\ln{x}}} \, dx\\
        &= \sum_{k=0}^{\infty} (-1)^k \int_0^1 \ln{(-\ln{x})}\frac{x^{2k}}{\sqrt{-\ln{x}}} \, dx.
    \end{align*}
    By Lemma 4,
    $$\int_0^1 \ln{(-\ln{x})}\frac{x^{2k}}{\sqrt{-\ln{x}}} \, dx = -[\gamma + 2\ln{2} + \ln{(2k+1)}]\sqrt{\frac{\pi}{2k+1}}.$$
    Thus,
    \begin{align}
        \int_0^1 \ln{(-\ln{x})}\frac{1}{\left(1+x^2\right)\sqrt{-\ln{x}}} \, dx &= \sqrt{\pi}\sum_{k=0}^{\infty} \frac{(-1)^{k+1}}{\sqrt{2k+1}}[\ln{(2k+1) + 2\ln{2} + \gamma}]\\
        &= \sqrt{\pi}\left(\beta'\left(\frac{1}{2}\right) - (2\ln{2} + \gamma)\beta\left(\frac{1}{2}\right)\right). \nonumber
    \end{align}
    Since $\beta(s)$ satisfies the functional equation
    $$\beta(1-s) = \left(\frac{2}{\pi}\right)^s\sin\left(\frac{\pi s}{2}\right)\Gamma(s)\beta(s)$$
    (see \cite[sec. 3.5]{functatlas}), differentiating both sides gives
    \begin{align*}
        -\beta'(1-s) = & \ln\left(\frac{2}{\pi}\right) \left(\frac{2}{\pi}\right)^s\sin\left(\frac{\pi s}{2}\right)\Gamma(s)\beta(s) + \frac{\pi}{2}\left(\frac{2}{\pi}\right)^s\cos\left(\frac{\pi s}{2}\right)\Gamma(s)\beta(s) \\
        &+ \left(\frac{2}{\pi}\right)^s\sin\left(\frac{\pi s}{2}\right)\Gamma(s)\psi(s)\beta(s) + \left(\frac{2}{\pi}\right)^s\sin\left(\frac{\pi s}{2}\right)\Gamma(s)\beta'(s).
    \end{align*}
    Let $s = \frac{1}{2}$,
    $$-\beta'\left(\frac{1}{2}\right) = \ln\left(\frac{2}{\pi}\right)\frac{\Gamma\left(\frac{1}{2}\right)}{\sqrt{\pi}}\beta\left(\frac{1}{2}\right) + \frac{\sqrt{\pi}}{2}\Gamma\left(\frac{1}{2}\right)\beta\left(\frac{1}{2}\right) + \frac{\Gamma\left(\frac{1}{2}\right)\psi\left(\frac{1}{2}\right)}{\sqrt{\pi}}\beta\left(\frac{1}{2}\right) + \frac{\Gamma\left(\frac{1}{2}\right)}{\sqrt{\pi}}\beta'\left(\frac{1}{2}\right).$$
    Since $\Gamma\left(\frac{1}{2}\right) = \sqrt{\pi}$ and $\psi\left(\frac{1}{2}\right) = -\gamma - 2\ln{2}$ (see \cite[sec. 5.4]{DLMF}),
    $$\beta'\left(\frac{1}{2}\right) = \left(-\frac{1}{2}\ln\left(\frac{2}{\pi}\right) - \frac{\pi}{4} + \frac{\gamma}{2} + \ln{2}\right)\beta\left(\frac{1}{2}\right).$$
    Substituting this into (8) gives
    $$\int_0^1 \ln{(-\ln{x})} \frac{1}{\left(1+x^2\right)\sqrt{-\ln{x}}}\,dx = -\left(\ln{\sqrt{\frac{8}{\pi}}} + \frac{\pi}{4} + \frac{\gamma}{2}\right)\sqrt{\pi}\,\beta\left(\frac{1}{2}\right).$$
\end{proof}

\section{Further Research}
As a generalization of the integrals evaluated in this paper, one could study the closed form of the classes of integrals of the forms
$$\int_0^1 \ln{(-\ln{x})}\,\frac{p(x)}{q(x)}\,dx$$
and
$$\int_0^1 \ln{(-\ln{x})}\, (-\ln{x})^s\,\frac{p(x)}{q(x)}\,dx,$$
where $p(x)$ and $q(x)$ are polynomials in terms of $x$. Also, similarly to the integral
$$\int_{\pi/4}^{\pi/2} \ln{\ln{\tan{x}}}\,dx,$$
which involves double logarithm and a trigonometric function inside, one could also try to evaluate the closed form of
$$\int_0^{\pi/2} \ln{\ln{\sec{x}}}\,dx.$$
More generally, one could study integrals of the form
$$\int_a^b \ln{\ln{f(x)}}\,dx,$$
where $f(x)$ is a function such that $f(x) >1$ on $(a,b)$.
\section{Disclosure Statement}
There are no competing interests to be declared.
\section{Data availability statement}
There is no data associated with this article.
\bibliographystyle{abbrv}
\bibliography{bibliography.bib}
\end{document}